\documentclass[review,3p,times]{elsarticle} 
\usepackage{setspace} 
\usepackage{subcaption}
\usepackage{graphicx}
\usepackage{epsfig}
\usepackage{epstopdf}
\usepackage{amssymb}
\usepackage{amsmath,bm}
\usepackage{units}
\usepackage{color}
\usepackage{multirow}
\usepackage{hyperref}
\usepackage{textgreek}
\usepackage{lineno}
\usepackage{float}
\usepackage{graphicx}
\usepackage[utf8]{inputenc}
\usepackage[labelfont=bf]{caption}
\usepackage{amssymb}
\usepackage{xurl}
\usepackage{booktabs}
\usepackage{mathtools,tikz,caption}
\DeclareRobustCommand\sampleline[1]{%
  \tikz\draw[#1] (0,0) (0,\the\dimexpr\fontdimen22\textfont2\relax)
  -- (2em,\the\dimexpr\fontdimen22\textfont2\relax);%
}

\biboptions{numbers,sort&compress}

\newcommand{\R}{\ensuremath{\mathbb{R}}}

\newcommand{\K}{\ensuremath{\mathbf{K}}}
\newcommand{\C}{\ensuremath{\mathbf{C}}}
\newcommand{\M}{\ensuremath{\mathbf{M}}}
\newcommand{\F}{\ensuremath{\mathbf{f}}}
\newcommand{\G}{\ensuremath{\mathbf{g}}}
\newcommand{\V}{\ensuremath{\mathbf{V}}}
\newcommand{\p}{\ensuremath{\mathbf{p}}}
\newcommand{\x}{\ensuremath{\mathbf{x}}}
\newcommand{\q}{\ensuremath{\mathbf{q}}}

\newcommand{\Krk}{\ensuremath{\mathbf{K}_{r,k}}}
\newcommand{\Crk}{\ensuremath{\mathbf{C}_{r,k}}}
\newcommand{\Mrk}{\ensuremath{\mathbf{M}_{r,k}}}
\newcommand{\Frk}{\ensuremath{\mathbf{f}_{r,k}}}
\newcommand{\Grk}{\ensuremath{\mathbf{g}_{r,k}}}

\newcommand{\pk}{\ensuremath{\mathbf{p}_{k}}}
\newcommand{\Vk}{\ensuremath{\mathbf{V}_{k}}}
\newcommand{\Tk}{\ensuremath{\mathbf{T}_{k}}}

\newcommand{\Krkt}{\ensuremath{\tilde{\mathbf{K}}_{r,k}}}
\newcommand{\Crkt}{\ensuremath{\tilde{\mathbf{C}}_{r,k}}}
\newcommand{\Mrkt}{\ensuremath{\tilde{\mathbf{M}}_{r,k}}}
\newcommand{\Frkt}{\ensuremath{\tilde{\mathbf{f}}_{r,k}}}
\newcommand{\Grkt}{\ensuremath{\tilde{\mathbf{g}}_{r,k}}}
\newcommand{\Vkt}{\ensuremath{\tilde{\mathbf{V}}_k}}

\journal{journal}
\begin{document}
\begin{frontmatter}
\title{Consistent Parametric Model Order Reduction by Matrix Interpolation for Varying Underlying Meshes}

\author[label1]{Sebastian Resch-Schopper \corref{cor1}}
\cortext[cor1]{Corresponding author at: sebastian.resch-schopper@tum.de}
\author[label2,label3]{Romain Rumpler}
\author[label1]{Gerhard Müller}

\address[label1]{Technical University of Munich, Chair of Structural Mechanics, Arcisstr. 21, 80333 Munich (sebastian.resch-schopper@tum.de, gerhard.mueller@tum.de)}
\address[label2]{KTH Royal Institute of Technology, MWL Laboratory for Sound and Vibration Research, Teknikringen 8, 10044 Stockholm (rumpler@kth.se)}
\address[label3]{KTH Royal Institute of Technology, Centre for ECO2 Vehicle Design and Digital Futures, Department of Engineering Mechanics, Teknikringen~8, 10044 Stockholm}

\begin{abstract}
Parametric model order reduction (pMOR) is a powerful tool for accelerating finite element (FE) simulations while maintaining parametric dependencies. For geometric parameters, pMOR by matrix interpolation is a well-suited approach because it does not require an affine representation of the parametric dependency, which is often not available for geometric parameters. However, the method requires that the underlying FE mesh has the same number of degrees of freedom and the same topology for all parameter configurations. This requirement can be difficult or even impossible to achieve for large parameter ranges or when automatic meshing is used. In this work, we propose a novel framework for pMOR by matrix interpolation for varying underlying meshes. The key idea is to understand the sampled reduced bases as continuous displacement fields that can be represented in different discretizations. By using mesh morphing and basis interpolation, the sampled reduced bases described in varying meshes can all be represented in terms of one reference mesh. This not only allows for performing pMOR by matrix interpolation, but also enables comparing the subspaces that the reduced bases span, which is important to detect strong changes that could lead to inconsistencies in the reduced operators. For mesh morphing, two strategies, namely morphing by spring analogy with elastic hardening and radial basis function morphing, were implemented and tested. Numerical experiments on a beam-shaped plate and a plate with a hole for one- and two-dimensional parameter spaces show that the proposed framework achieves high accuracy for both morphing methods and performs significantly better than two existing approaches for pMOR by matrix interpolation for varying underlying meshes. 
\end{abstract}


\end{frontmatter}


\section{Introduction}%
\label{sec:Introduction}

The Finite Element Method (FEM) has become a widely used tool to simulate the behavior of technical systems. However, for dynamic problems, these simulations can become computationally very expensive because large systems of equations must be solved repeatedly for various instances in time or frequency. For parametric applications, the computational complexity increases further because the systems must be evaluated for different configurations of the parameters. To reduce the computational effort of FE simulations, projection-based model order reduction can be used. In these methods, the system is projected onto a lower-dimensional subspace described by a reduced basis, in which the full solution can be approximated well and solved with significantly less computational effort. \cite{Benner2015} \\
To incorporate parametric dependencies in these reduced models, parametric model order reduction (pMOR) can be used. These methods can be distinguished between global and local methods \cite{Benner2015}. In the former, one reduced basis is computed for the complete parameter space, which is used to project the full system. The latter methods instead interpolate reduced quantities such as the reduced bases \cite{Amsallem2008, Goutaudier2023}, the reduced operators \cite{Amsallem2009, Panzer2010, Amsallem2011}, or the reduced transfer function \cite{Ionita2014, Feng2022, Rodriguez2023}. For global methods to be efficient, an affine representation of the parametric dependency of the full operators, i.e., a representation in terms of (a few) constant matrices multiplied by scalar functions, must be known, so that the projection does not have to be performed in the online phase. This representation can often be easily achieved for material parameters, but might not be available for geometric parameters. \\
In this work, we focus on geometric parameters and assume that no affine decomposition of the parametric dependency is available. Therefore, we concentrate on parametric model order reduction by matrix interpolation \cite{Panzer2010}, which does not require an affine representation of the parametric dependency. In this method, the reduced operators are interpolated, which leads to a small computational effort in the online phase, as the projection onto the reduced basis does not have to be performed. Furthermore, the dependency of the reduced operators might be approximated well with simple functions like polynomials, contrary to the dependency of the transfer function, which is a rational function \cite{Geradin2015}. \\
In pMOR by matrix interpolation, reduced systems are first individually computed for a set of sample points. Eventually, the entries of the reduced operators shall be interpolated. In order to do so, a prior transformation to a common coordinate system is necessary so that the interpolation of the entries is meaningful. This transformation, however, is only possible if the collected reduced bases are similar for all samples so that the basis vectors can be rearranged via a linear transformation to resemble a common coordinate system. If the subspaces that the individually sampled reduced bases span differ strongly, however, individual basis vectors cannot be matched. This leads to inconsistencies in the reduced operators, which prevent a meaningful interpolation of the entries and lead to high errors in the predicted systems. Several causes and remedies to reduce these inconsistencies are discussed in \cite{Resch-Schopper2024}. In that work, we proposed an adaptive sampling algorithm and a partitioning of the parameter space into separate regions so that the collected reduced operators are consistent with their neighboring samples in this region. \\
A general challenge when using geometric parameters is the treatment of the FE mesh. All pMOR methods, except the interpolation of the reduced transfer function, require that the size of the full system stays the same. In pMOR by matrix interpolation, for example, this is required in order to compute a common basis and compute the transformation matrix that is applied to the sampled reduced operators. Maintaining the number of nodes and topology of the mesh when using geometric parameters can be achieved by utilizing a reference mesh that is adjusted to the queried geometric configuration through mesh morphing. Overviews of these methods can be found in \cite{Alexa02, Selim2016}, and applications in the context of pMOR in \cite{Burgard13, Manzoni2012, Manzoni2017, Lupini2019, Agathos2020, Shengfang2020, Mencik2021, Mencik2024}. For large variations of the geometry, however, morphing the mesh can lead to highly distorted elements, which introduce errors in the FE simulations \cite{Agathos2020}. \\
To circumvent the problem of distorted meshes, two approaches for using pMOR by matrix interpolation for varying underlying meshes have been proposed in the literature. In \cite{Amsallem2016}, it is suggested to perform the transformation based on the reduced operators instead of the reduced bases. This is done by computing the transformation matrix via an optimization problem in which the difference between a reference reduced operator and the sampled, transformed reduced operator is minimized. In \cite{Roy2021}, an alternative approach is suggested, in which zero matrices are appended to the sampled reduced bases so that they are all the same size, which allows performing the required operations. Both approaches, however, are not able to detect inconsistencies in the reduced bases. The approach proposed in \cite{Amsallem2016} is only based on the reduced operators and thus does not take the reduced bases into account. In the approach suggested in \cite{Roy2021}, zero-padding the reduced bases influences them in a way that inconsistencies in them cannot be detected. \\
We therefore propose a novel method for pMOR by matrix interpolation for varying underlying meshes that enables the detection of inconsistencies in the reduced bases. The key point of this approach is to understand the reduced basis vectors as continuous displacement fields that are represented in a discrete way via the underlying mesh. To describe them all with respect to one mesh, we define a reference mesh, in which we want to represent the continuous displacement fields described by all sampled reduced basis vectors. This is achieved by morphing the reference mesh to the geometries of the sampled meshes and then evaluating the sampled reduced bases at the nodes of the morphed reference mesh. This way, all operations required in pMOR by matrix interpolation can be performed because the sampled reduced bases are all described in terms of the reference mesh. The same principle of performing mesh morphing and basis interpolation has already been applied for learning solutions of PDEs using neural networks \cite{Casenave2023}, but not for pMOR by matrix interpolation so far, to our knowledge. Furthermore, we give comments on some aspects regarding the usage of the approach in the context of pMOR by matrix interpolation. \\
The remainder of this paper is structured as follows: In section \ref{sec:Background}, the mathematical description of the investigated systems and projection-based model order reduction are introduced. Parametric model order reduction by matrix interpolation for common underlying meshes is explained in section \ref{sec:pMOR_MI}. In section \ref{sec:pMOR_varyingMeshes}, the proposed approach for varying underlying meshes is presented. Finally, the results of the proposed approach are presented in section \ref{sec:Results}, and a conclusion is given in section \ref{sec:conclusion}.

\section{Theoretical Background}
\label{sec:Background}

\subsection{Problem Definition}
\label{sec:ProblemDefinition}

We consider parameter-dependent, linear time-invariant systems in second-order form with single input and single output obtained from a discretization with the finite element method (FEM). By applying a Laplace transformation, the system is obtained in the frequency domain:
\begin{equation}
   \Sigma(\p): \begin{cases}
		s^2\M(\p) \q(\p,s) + s\C(\p) \q(\p,s) + \K(\p) \q(\p,s) &= \F(\p) u(s), \\
		\hfill y(\p,s) &= \G(\p)\q(\p,s).
		\end{cases} \label{eq:FOM}
\end{equation}
The system is assumed to depend on $d$ parameters $\p = [p_1, p_2, \dots, p_d] \in \Omega \subset \R^d$, where $\Omega$ denotes a bounded domain. Furthermore, we assume that the discretization of the finite element (FE) model may change with the parameters so that the number of degrees of freedom $n$ depends on $\p$. The matrices $\M$, $\C$, $\K$: $\Omega \rightarrow \mathbb{R}^{n(\p) \, \times \, n(\p)}$ denote the parameter-dependent mass, damping, and stiffness matrix, $\q \in \mathbb{R}^{n(\p)}$ the degrees of freedom, $s \in \mathbb{C}$ the complex frequency, $u(s) \in \mathbb{R}$ and $\F$ : $\Omega \rightarrow \mathbb{R}^{n(\p) \, \times \, 1}$ the input and $y(\mathbf{p},s) \in \mathbb{R}$ and $\G$ : $\Omega \rightarrow \mathbb{R}^{1 \, \times \, n(\p)}$ the output.

\subsection{Projection-based Model Order Reduction}
\label{sec:MOR}

For complex problems, the system in Eq.~\eqref{eq:FOM} can comprise thousands or millions of degrees of freedom and thus become computationally very expensive to solve for multiple instances in frequency. In order to reduce the computational effort for a specific parameter point $\p_k$, projection-based model order reduction (MOR) may be used. The idea of MOR is to project the system of Eq.~\eqref{eq:FOM} onto a lower-dimensional subspace represented by the reduced basis $\V_k \in \mathbb{R}^{n(\p_k) \times r}$. This leads to the following reduced system of size $r \ll n(\p_k)$:

\begin{equation}
   \Sigma_r(\p_k): \begin{cases}
		s^2\M_{r,k} \q_r(s) + s\C_{r,k} \q_r(s) + \K_{r,k} \q_r(s) &= \F_{r,k} u(s), \\
		\hfill y_r(s) &= \G_{r,k} \q_r(s),
		\end{cases} \label{eq:ROM}
\end{equation}
where
\begin{equation}
\begin{aligned}
    \Mrk &= \Vk^\intercal \M(\pk) \Vk, \qquad &\Crk &= \Vk^\intercal \C(\pk) \Vk, \qquad \Krk = \Vk^\intercal \K(\pk) \Vk \\
    \Frk &= \Vk^\intercal \F(\pk), \qquad &\Grk &= \G(\pk) \Vk.
    \label{eq:Projection}
\end{aligned}
\end{equation}
For computing the reduced basis $\V_k$, numerous methods have been developed over the past decades, such as modal methods \cite{Tiso2021}, moment matching \cite{Benner2021}, or Proper Orthogonal Decomposition (POD) \cite{Antoulas2005}. In the context of pMOR by matrix interpolation, modal methods have beneficial properties \cite{Fischer2015, Resch-Schopper2024} and have been advantageously used in several contributions on pMOR by matrix interpolation \cite{Panzer2010, Fischer2014, Fischer2015, Mencik2021, Mencik2024}. We thus focus on modal truncation only in the following. \\
For computing the reduced basis with modal truncation, the undamped eigenvalue problem
\begin{equation}
    (\K_k - \omega^2 \M_k) \pmb{\Phi}_k = 0
\end{equation}
must be solved, where $\omega$ denotes the angular eigenfrequencies of the system and $\pmb{\Phi}_k$ are the associated eigenvectors~\cite{Tiso2021}. Often, all eigenvectors within twice the frequency range of interest are used as the reduced basis. Alternatively, the $r$ eigenmodes that contribute most to the input-output-behavior based on a dominance criterion may be selected \cite{Rommes2008}.

\section{Parametric Model Order Reduction by Matrix Interpolation for Matching Underlying Meshes}%
\label{sec:pMOR_MI}

PMOR approaches aim at preserving the parametric dependency of the full-order model (FOM) $\Sigma(\p)$ in the reduced-order model (ROM). Many pMOR approaches require that an (efficient) affine decomposition of the parametric dependency is available. In this work, however, we assume that there is no such representation available. Therefore, we use pMOR by matrix interpolation \cite{Panzer2010}. We will first revise this method for the case where the topology of the mesh and the number of degrees of freedom do not change with the parameters. \\
The first step of pMOR by matrix interpolation is to sample the FOM for $K$ sample points $\p_1, \dots, \p_K$. Each of the sampled FOMs is then reduced individually by projection-based model order reduction. It is important to note that in the original version of pMOR by matrix interpolation, the reduced size $r$ of the local ROMs must be the same for all samples. An extension of this approach to interpolate among differently sized reduced models has been presented in \cite{Geuss2014}. This approach, however, leads to either a deteriorated accuracy or singular transformed reduced matrices for some sample points , which is why it is not further regarded here. In the end, the entries of the reduced operators shall be interpolated so that the reduced system can be predicted for queried parameter points. At this stage, however, the sampled ROMs might be described in different reduced bases, because the reduction was performed for each FOM individually. Therefore, the local ROMs must be transformed to a common coordinate system first. To compute this common coordinate system, a singular value decomposition (SVD) is performed on the library of all collected reduced bases:
\begin{equation}
    \mathbf{U}\pmb{\Sigma}\mathbf{Y} = [\V_1, \V_2, \dots, \V_K]. \label{eq:Vg}
\end{equation}
The common coordinate system $\mathbf{R}$ is then obtained as the first $r$ columns of $\mathbf{U}$. Subsequently, the local reduced systems shall be transformed so that the reduced bases they are described in resemble $\mathbf{R}$ as much as possible. For that, a transformation matrix defined as
\begin{equation}
    \Tk = (\mathbf{R}^\intercal \Vk)^{-1} \label{eq:Tk}
\end{equation}
is applied to the reduced operators:
\begin{equation}
\begin{aligned}
    \Mrkt &= \Tk^\intercal \Mrk \Tk, \qquad &\Crkt &= \Tk^\intercal \Crk \Tk, \qquad \Krkt = \Tk^\intercal \Krk \Tk    \\
    \Frkt &= \Tk^\intercal \Frk, \qquad &\Grkt &= \Grk \Tk.
    \label{eq:TransfOp}
\end{aligned}
\end{equation}

Finally, an interpolation or regression model for the reduced operators can be generated either for each entry separately or for the complete operators. To ensure that during the interpolation, important properties such as positive definiteness are preserved, the interpolation can be carried out on the tangent space of the manifold of symmetric positive definite matrices \cite{Amsallem2011}. Alternatively, the Cholesky factors of the reduced operators can be interpolated  \cite{Amsallem2016}.

\subsection{Inconsistencies in the Reduced Bases}

After the transformation, each local ROM is described in the transformed reduced basis $\Vkt$, which is given by $\Vkt = \Vk \Tk$. We can thus write the complete parametric dependency to be interpolated as
\begin{equation}
    \tilde{\mathbf{K}}_r(\p) = \tilde{\V}(\p)^\intercal \K(\p) \tilde{\V}(\p). \label{eq:Krt(p)}
\end{equation}
If the sampled local reduced bases $\Vk$ are similar, a reordering of the basis vectors in the course of the transformation is effective to make the interpolation of the reduced operators meaningful. However, this is not the case if the subspaces spanned by the sampled reduced bases differ greatly, because the linear transformation does not change the subspace. Consequently, strong changes in the reduced bases, which will be termed inconsistencies in the following, may greatly affect the total parametric dependency of the transformed reduced operators. In our previous work \cite{Resch-Schopper2024}, we discuss possible causes for inconsistencies, such as the choice of the MOR method or the phenomenon of mode switching and truncation, in detail. \\
In order to detect such inconsistencies, the angle between the subspaces spanned by the reduced bases can be measured. For two orthonormal bases $\mathbf{V}_i$ and $\mathbf{V}_j$, this is done by first performing the following SVD \cite{Amsallem2016}:
\begin{equation}
    \mathbf{W} \pmb{\Sigma} \mathbf{Z} = \V_i^\intercal \V_j,
    \label{eq:SubspaceAngles}
\end{equation}
The angles between the subspaces $\theta_l$ are then obtained from the singular values in $\pmb{\Sigma} = \operatorname{diag}(\sigma_1, \dots, \sigma_r)$ \cite{Amsallem2016}:
\begin{equation}
    \theta_l = \arccos(\sigma_l), \quad l = 1, \dots r.
\end{equation}
Large values for the angles indicate strong differences in the spanned subspaces. \\
For removing these inconsistencies in the parametric ROM, we proposed an adaptive sampling and clustering algorithm in \cite{Resch-Schopper2024}. The algorithm starts with a small set of initial samples and subsequently adds more samples until the subspace angles between the reduced bases of neighboring samples are either below a lower threshold, indicating consistency, or above an upper threshold, indicating inconsistency. By this, clusters of samples can be identified, so that neighboring samples within a cluster are consistent with each other. Finally, the steps of pMOR by matrix interpolation can be performed for each cluster separately leading to multiple local pROMs.

\section{Parametric Model Order Reduction for Varying Underlying Meshes}%
\label{sec:pMOR_varyingMeshes}

\begin{figure}
    \centering
    \includegraphics[width=0.5\linewidth]{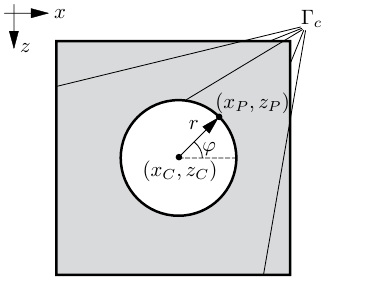}
    \caption{Visualization of the characteristic features for a plate with a circular hole.}
    \label{fig:CharacteristicPoints}
\end{figure}

A prerequisite of the original version of pMOR by matrix interpolation is that the number of degrees of freedom of the FOM must stay the same for all parameter configurations. Otherwise, the concatenation of the sampled local reduced bases in Eq.~\eqref{eq:Vg} and the computation of the transformation matrix in Eq.~\eqref{eq:Tk} are not possible. Furthermore, it is implicitly assumed that the position of the nodes relative to the geometry stay the same, i.e. a node in the center of the structure will stay in the center when the parameters change. This is necessary so that the comparison of the sampled local reduced bases is meaningful. Finally, the angle between the subspaces spanned by the reduced bases can only be computed when the sampled reduced bases are described in the same high-dimensional space. This is important for detecting inconsistencies stemming from, e.g., mode switching and truncation. However, for several applications, these prerequisites might not be fulfilled. When using geometric parameters, for example, using the same mesh for each parameter configuration can lead to extremely distorted meshes, which in turn leads to inaccurate results in the FE simulations \cite{Agathos2020}. A further example is when automatic meshing is used. In these cases, there might be no possibility to ensure that the number of degrees of freedom and the relative position of the nodes stay the same. We thus want to extend the approach of pMOR by matrix interpolation for cases where the underlying FE mesh may vary with the parameters. We assume that both the number of degrees of freedom and the relative position of the nodes within the structure may change with the parameters. However, the meshes are not allowed to be arbitrary. Instead, we put the mild restriction that there must be specific features in the structure that are present for each parameter configuration, such as boundaries and corners, which we call characteristic points of the geometry in the following. The remaining parts of the structure will be referred to as non-characteristic points. This is visualized in Figure \ref{fig:CharacteristicPoints} for a plate with a hole, for which the radius of the hole is used as a parameter. Independent of the size of the radius, the plate will comprise the four bounding edges and the boundary of the circular hole, which thus make up the characteristic features of the structure $\Gamma_c$. Furthermore, it must be known analytically how the coordinates of the characteristic points change with a change of the parameters. In the  plate shown in Figure \ref{fig:CharacteristicPoints}, for example, a point on the circumference of the hole with coordinates $x_P$ and $z_P$ is parametrized by
\begin{align}
    x_P &= x_C + r \cdot \cos(\varphi), \\
    z_P &= z_C + r \cdot \sin(\varphi),
\end{align}
where $x_C$ and $z_C$ are the coordinates of the center of the hole and $\varphi = \arctan2(x_P - x_C, z_P - z_C)$. 

\subsection{Previous work}

The problem of pMOR by matrix interpolation for arbitrary underlying meshes was first treated in \cite{Amsallem2016}. In this work, the transformation matrix $\mathbf{T}_k$ is computed based on the reduced operators instead of the reduced bases as shown in Eq.~\eqref{eq:Tk}. This is done by solving an optimization problem such that the difference between the transformed reduced operators and one of the sampled reduced operators, which is chosen as a reference, is minimized. According to the total parametric dependency of the transformed reduced operators shown in Eq.~\eqref{eq:Krt(p)}, there are two factors that enter the optimization criterion, which might conflict with each other: the parametric dependency of the reduced bases and the parametric dependency of the reduced system matrices. Only the former should be minimized, while the latter is inherent to the system and should be maintained. \\
An alternative and computationally less expensive approach was developed in \cite{Roy2021}. To be able to perform the concatenation of the reduced bases in Eq.~\eqref{eq:Vg}, zero matrices are appended to the sampled reduced bases so that all reduced bases have the same size. \\
Both approaches, however, are not suited to detect inconsistencies in the reduced bases when the subspaces they span vary greatly. This is because the approach presented in \cite{Amsallem2016} does not take the reduced bases into account and therefore cannot distinguish between changes in the reduced bases and changes in the reduced operators. In the approach proposed in \cite{Roy2021}, the measured subspace angles would not be meaningful because padding the reduced bases with zeros influences the subspace they span. Therefore, we propose an alternative approach for pMOR by matrix interpolation for varying meshes that takes the reduced bases into account without affecting them.

\subsection{Proposed approach}

\begin{figure*}[htb]
	\centering
	\begin{subfigure}[t]{0.3\textwidth}
		\centering
		\includegraphics[width=\textwidth]{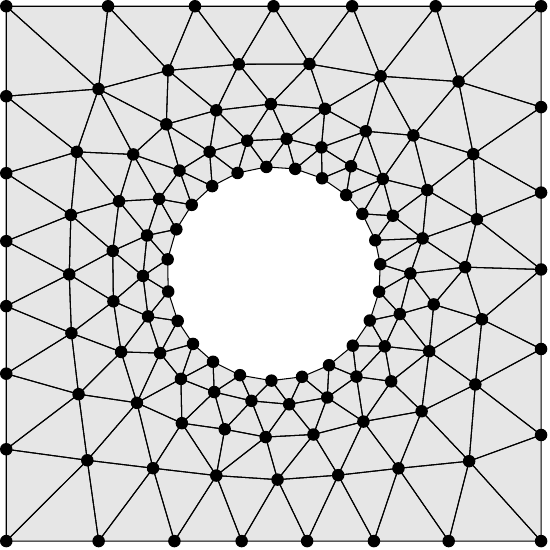}
		\caption{Mesh of one of the samples.}  
		\label{fig:SampledMesh}
	\end{subfigure}
	\hfill
	\begin{subfigure}[t]{0.3\textwidth}
		\centering 
		\includegraphics[width=\textwidth]{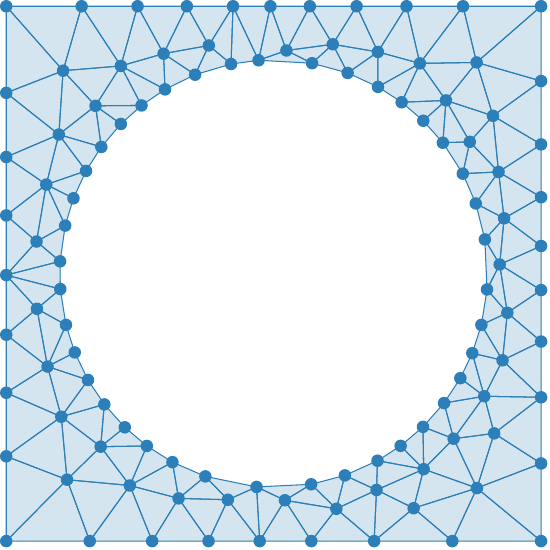}
		\caption{Reference mesh.}  
		\label{fig:ReferenceMesh}
	\end{subfigure}
    \hfill
    \begin{subfigure}[t]{0.3\textwidth}
		\centering 
		\includegraphics[width=\textwidth]{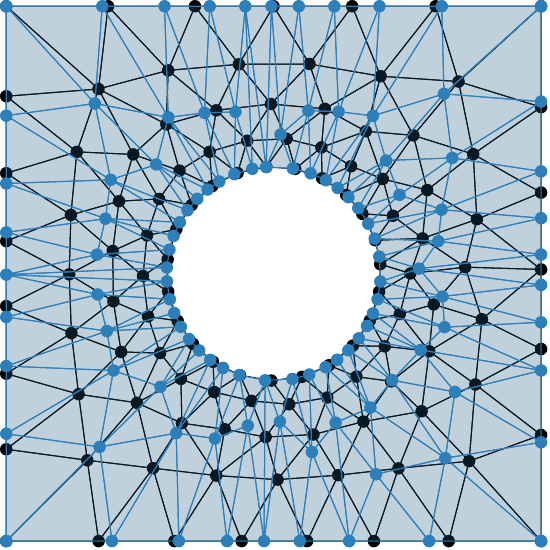}
		\caption{Reference mesh morphed to the geometry of the sampled mesh so that the characteristic parts of the two meshes overlap.}  
		\label{fig:MorphedMesh}
	\end{subfigure}
	\caption{Visualization of mesh morphing.}  
	\label{fig:MeshMorphing}
\end{figure*}

\begin{figure*}[htb]
	\centering
    \hspace{2cm}
	\begin{subfigure}[t]{0.3\textwidth}
		\centering
		\includegraphics[width=\textwidth]{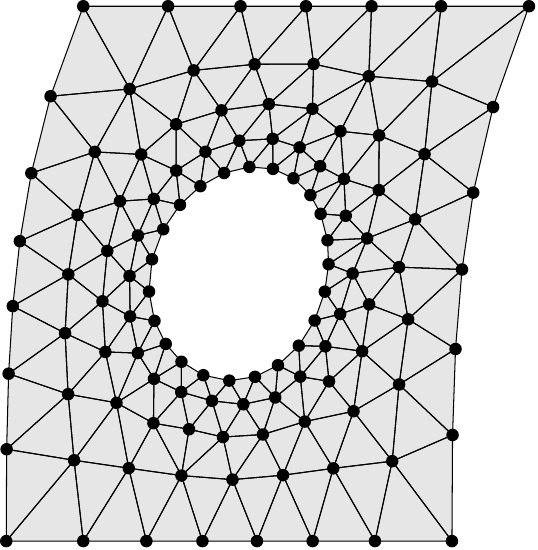}
		\caption{Visualization of one basis vector represented in terms of the sampled mesh.}  
		\label{fig:DisplacedSampledMesh}
	\end{subfigure}
	\hfill
	\begin{subfigure}[t]{0.3\textwidth}
		\centering 
		\includegraphics[width=\textwidth]{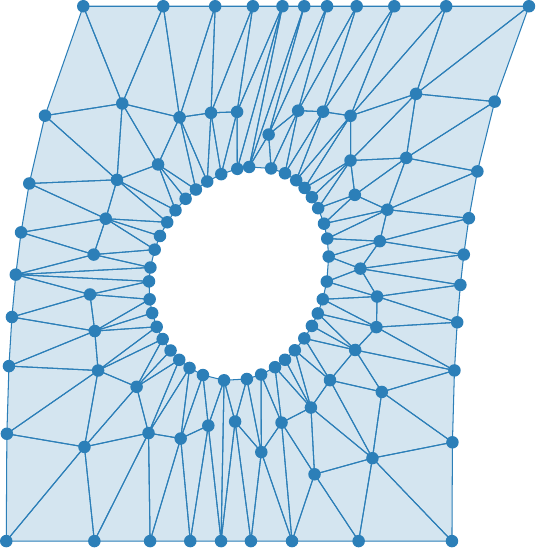}
		\caption{Basis vector of sampled mesh represented in terms of the morphed reference mesh.}  
		\label{fig:DisplacedReferenceMesh}
	\end{subfigure}
    \hspace{2cm}
    \caption{Visualization of representing a basis vector described in one of the sampled meshes in terms of the reference mesh.}  
	\label{fig:DisplacementFieldMeshes}
\end{figure*}

The key idea of our proposed approach is to understand the vectors of the reduced basis as continuous displacement fields, which are represented in a discrete form due to the underlying mesh. With the help of the shape functions used in the FE discretization, this continuous displacement field can be retrieved from the discrete representation. In order to make the displacement fields described by the sampled local reduced bases comparable, they should all be represented on the same mesh. Therefore, we first define the mesh of one of the sampled parameter configurations as the reference mesh. This can be chosen arbitrarily from all sampled geometries, but we will later give some considerations regarding this choice. After having computed the local reduced bases for all samples, they shall all be represented in terms of the reference mesh. For that, the reference mesh is morphed for each individual sampled mesh such that the characteristic features of the morphed mesh overlap with the characteristic features of the sampled mesh under consideration. This is illustrated in Fig. \ref{fig:MeshMorphing} for the plate with a hole shown in Fig. \ref{fig:CharacteristicPoints}. Standard mesh morphing techniques can be used in this step, which will be discussed in detail in the next section. When the characteristic features of the morphed and the sampled mesh overlap, the continuous displacement fields represented by the reduced bases in the sampled mesh can be evaluated at the nodes of the morphed reference mesh. This step will be explained in further detail in section \ref{sec:BasisInterpolation}. Fig. \ref{fig:DisplacedReferenceMesh} illustrates a displacement field described by a reduced basis vector represented in terms of the sampled and the morphed reference mesh. Finally, all sampled reduced bases are described in the reference mesh, so that a concatenation and comparison of the bases is possible. In the following, we will explain both mesh morphing and basis interpolation in more detail. \\

\subsection{Mesh Morphing}

The purpose of mesh morphing is to transform the nodes of a sampled mesh to the geometry of the reference mesh such that the characteristic features of the morphed and the reference mesh overlap. For this task, a variety of approaches have been developed, such as for example morphing based on linear elasticity \cite{Tezduyar1992,Nielsen2002,Yang2007}, spring analogies \cite{Batina1990,Farhat1998}, barycentric mapping \cite{Tutte1963} and data-driven morphing \cite{deBoer2007,Witteveen2010}. In the following, we focus on morphing by the spring analogy and radial basis function morphing, which is a form of data-driven morphing. In preliminary studies conducted in the course of this work, these two approaches proved to work most robustly, also for large changes in the geometries. Furthermore, we assume that the mesh consists of triangular elements, which are most commonly used in automatic meshing. 

\subsubsection{Spring Analogy with Elastic Hardening}

\begin{figure}
    \centering
    \includegraphics[width=0.5\linewidth]{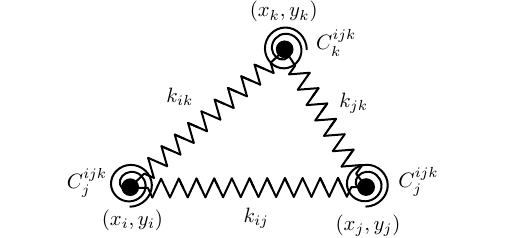}
    \caption{Triangular element with lineal and torsional springs.}
    \label{fig:ElementSpringAnalogy}
\end{figure}

Mesh morphing can be performed by solving a problem of linear elasticity: The displacement $\bar{\q}$ of nodes on the characteristic parts of the mesh can be imposed as Dirichlet boundary conditions and the new position of the remaining nodes can be computed by solving a static problem: 
\begin{align}
    \K \q &= 0, \\
    \q &= \bar{\q} \quad \text{on } \Gamma_c. 
\end{align}
The imposed displacement is thereby incrementally applied in multiple steps and the stiffness matrix is updated after each increment \cite{Nielsen2002}. The number of steps is problem-dependent but should be chosen large enough to avoid collapsing elements. For complex structures, reassembling the stiffness matrix multiple times can be computationally very expensive. Furthermore, this approach might lead to collapsing elements if the volume of the mesh is substantially affected in each step. In order to prevent both problems, an approach based on a spring analogy has been presented in \cite{Farhat1998}. In this approach, the edges of all elements are replaced by lineal springs as shown in Fig. \ref{fig:ElementSpringAnalogy}, whose stiffness $k_{ij}$ is inversely proportional to the length $l_{ij}$ of the edge between the nodes $i$ and $j$ with coordinates $(x_i \, | \, y_i)$ and $(x_j \, | \, y_j)$: 
\begin{equation}
    k_{ij} \propto \frac{1}{l_{ij}}.
\end{equation}
The stiffness of such a spring element in 2D is thus given by
\begin{equation}
    \mathbf{K}_{lineal}^{ij} = \frac{1}{l_{ij}} \begin{bmatrix} \cos^2(\alpha) & \sin(\alpha) \cos(\alpha) & -\cos^2(\alpha) & -\sin(\alpha) \cos(\alpha) \\
    \sin(\alpha) \cos(\alpha) & \sin^2(\alpha) & -\sin(\alpha) \cos(\alpha) & -\sin^2(\alpha) \\
    -\cos^2(\alpha) & -\sin(\alpha) \cos(\alpha) & \cos^2(\alpha) & \sin(\alpha) \cos(\alpha) \\
    -\sin(\alpha) \cos(\alpha) & -\sin^2(\alpha) & \sin(\alpha) \cos(\alpha) & \sin^2(\alpha) \end{bmatrix},
\end{equation}
where $\alpha = \tan^{-1}\left(\frac{y_j - y_i}{x_j - x_i}\right)$. As the imposed deformation reduces $l_{ij}$, the spring stiffness increases, preventing one node from being pushed through another. Nevertheless, a node could be pushed through an opposite edge. In order to prevent this, torsional springs are added to the nodes, whose stiffness $C_i^{ijk}$ is computed based on the angle between two neighboring edges: 
\begin{equation}
    C_i^{ijk} = \frac{l_{ij}^2 l_{ik}^2}{4A_{ijk}^2},
\end{equation}
where $l_{ij}$ and $l_{ik}$ are the lengths of the edges between nodes $i$ and $j$ respectively $i$ and $k$, and $A_{ijk}$ is the area of the triangle formed by nodes \( i \), \( j \), and \( k \). The stiffness matrix of a triangular element induced by the torsional springs is then computed by
\begin{equation}
    \mathbf{K}_{torsional}^{ijk} = \mathbf{R}^{{ijk}^\intercal} \mathbf{C}^{ijk} \mathbf{R}^{ijk},
\end{equation}
where 
\begin{equation}
    \mathbf{C}^{ijk} = \begin{bmatrix}
        C_i^{ijk} & 0 & 0 \\ 0 & C_j^{ijk} & 0 \\ 0 & 0 & C_k^{ijk}
    \end{bmatrix},
\end{equation}
and 
\begin{equation}
    \mathbf{R}^{ijk} = \begin{bmatrix}
        b_{ik} - b_{ij} & a_{ij} - a_{ik} & b_{ij} & - a_{ij} & - b_{ik} & a_{ik} \\
        -b_{ji} & a_{ji} & b_{ji} - b_{jk} & a_{jk} - a_{ji} & b_{jk} & - a_{jk} \\
        b_{ki} & - a_{ki} & -b_{kj} & a_{kj} & b_{kj} - b_{ki} & a_{ki} - a_{kj}
    \end{bmatrix},
\end{equation}
with $a_{ik} = \frac{x_k - x_i}{l_{ik}^2}$ and $b_{ik} = \frac{y_k - y_i}{l_{ik}^2}$. \\
To further prevent these elements from collapsing during the morphing process, we follow an empirical approach proposed in \cite{Moavenian2016,Nazem2021}. There, an elastic hardening coefficient $c_{EH}$ is computed and used to scale the torsional springs. This coefficient is computed based on the circumcircle $R$ and the incircle $r$ of each triangle, which are given by the following formulas
\begin{align}
    R &= \frac{l_{ij} l_{jk} l_{ik}}{\sqrt{(l_{ij}+l_{jk}+l_{ik})(l_{ij}+l_{jk}-l_{ik})(l_{ij}-l_{jk}+l_{ik})(-l_{ij}+l_{jk}+l_{ik})}}, \\[2mm]
    r & = \frac{\sqrt{a(a-l_{ij})(a-l_{jk})(a-l_{ik})}}{a},
\end{align}
where $a = \frac{1}{2}(l_{ij} + l_{jk} + l_{ik})$. In \cite{Moavenian2016}, several empirical functions are proposed for computing $c_{EH}$; in this work, we are using $c_{EH} = 4 \frac{R}{r} - 1$, which is then used to scale the torsional stiffness matrix $\mathbf{K}_{torsional}^{ijk}$. The global stiffness matrix based on this spring analogy is obtained by assembling the edge-based and elemental stiffness matrices $\mathbf{K}_{lineal}^{ij}$ and $\mathbf{K}_{torsional}^{ijk}$.

\subsubsection{Radial Basis Function Morphing}

Radial basis functions are a method for interpolating scattered data by a weighted sum of radially symmetric kernel functions \cite{Hardy1971}. In the context of mesh morphing, this can be used to predict the position of the non-characteristic nodes of the morphed mesh based on the known prescribed displacement of the characteristic nodes \cite{deBoer2007}. More specifically, we want to predict the displacement $\q$, which the non-characteristic nodes have to undergo in the morphing. The radial basis function interpolant takes the form 
\begin{equation}
    q(\mathbf{x}) = \sum_{i=1}^{N} \gamma_i \psi(\| \mathbf{x} - \mathbf{x}_i \|_2) + f(\x),
\end{equation}
where $N$ is the number of nodes on $\Gamma_c$, $\gamma_i$ are the weights of the kernel functions $\psi$, $\x$ are the coordinates of an arbitrary node, $\x_i$ are the coordinates of the nodes on $\Gamma_c$, $\| \cdot \|_2$ denotes the Euclidian norm and $f$ is a polynomial. For the kernel function, any radially symmetric function can be used. Some of those use a shape parameter that must be tuned to the specific problem at hand, which can be a difficult task. Therefore, we focus on polyharmonic kernel functions in this work, which do not use a shape parameter. For a given order $m$, this kernel function is defined as
\begin{equation}
    \psi(r) = \begin{cases} r^m \quad &\text{ for } m = 1, 3, 5, \dots \\ 
    r^m \ln(r) \quad &\text{ for } m = 2,4,6, \dots \end{cases}
\end{equation}
Furthermore, we are using linear polynomials to ensure that rigid body modes are covered. The interpolant then writes
\begin{equation}
    q(\mathbf{x}) = \sum_{i=1}^{N} \gamma_i \psi(\| \mathbf{x} - \mathbf{x}_i \|_2) + \mathbf{w} \begin{bmatrix} 1 \\ \x \end{bmatrix},
\end{equation}
where $\mathbf{w}$ are the weights of the polynomial terms. This specific form corresponds to a polyharmonic spline. \\\
To obtain the weights $\pmb{\gamma}$ and $\mathbf{w}$, the following system of equations must be solved:
\begin{equation}
    \begin{bmatrix}
        \mathbf{A} & \mathbf{B} \\ \mathbf{B}^\intercal & \mathbf{0} 
    \end{bmatrix} \begin{bmatrix}
        \pmb{\gamma} \\ \mathbf{w}
    \end{bmatrix} = \begin{bmatrix}
        \bar{\mathbf{q}} \\ \mathbf{0} 
    \end{bmatrix}, \label{eq:RBF_ES}
\end{equation}
where 
\begin{equation}
    A_{ij}=\psi(\|\x_i - \x_j\|_2), \quad \mathbf{B} = \begin{bmatrix} 1 & 1 & \dots & 1 \\ \x_1^\intercal & \x_2^\intercal & \dots & \x_N^\intercal \end{bmatrix}^\intercal
\end{equation}
In the context of mesh morphing, the weights are fitted for each component of the displacement $\mathbf{q}$ separately. Afterwards, the displacement components of the non-characteristic nodes are predicted using the fitted interpolant and the position of the morphed nodes is obtained by $\x + \mathbf{q}(\x)$.

\subsection{Basis Interpolation} 
\label{sec:BasisInterpolation}

When the characteristic points of the morphed and the reference mesh overlap, the displacement field represented by the basis vectors must be evaluated at the nodes of the other mesh. For that, the shape functions can be used as the displacement at any point within an element is given by 
\begin{equation}
    \q(\x) = \sum_{i=1}^{n_e} N_i(\x) \q_i,
\end{equation}
where $n_e$ is the number of nodes of the element, $N_i$ are the shape functions and $\q_i$ the vector of nodal displacements. Usually, the shape functions are only available in terms of the natural coordinates $\pmb{\xi}$ of the elements. Therefore, the physical coordinates must first be transformed into natural coordinates via the relation
\begin{equation}
    \x = \sum_{i=1}^{n_e} N_i(\pmb{\xi}) \x_{e,i}, \label{eq:CoordinateTransformation}
\end{equation}
where $\x_{e,i}$ correspond to the coordinates of the nodes of the element. For triangular elements with linear shape functions, there is a unique solution for the inverse relationship of Eq.~\eqref{eq:CoordinateTransformation}; for elements with higher-order shape functions, a non-linear system of equations must be solved to obtain the natural coordinates from the physical coordinates.

\subsection{Implementational Details}
\label{sec:ImplementationalDetails}

In this section, we provide details on the implementation of the methods described in the previous sections. Whenever operations on two or more reduced bases that are represented in different meshes shall be performed, a reference mesh must be chosen first. We propose always choosing the mesh with the highest number of nodes as the reference mesh because the reduced bases of the other meshes will be evaluated in this mesh. If a rather coarse mesh were chosen, it is possible that high-frequency information of the reduced bases of denser meshes could be lost. This is avoided if the densest mesh is chosen as reference. \\
Regarding the mesh morphing, either the sampled meshes can be morphed to the shape of the reference mesh, or vice versa. After the morphing, the reduced basis will be evaluated at the nodes of the other mesh by using the shape functions. In order to avoid distorted meshes in this step, the reference mesh should be morphed to the shape of the sampled meshes, in which the reduced bases are given. \\
Generally, for evaluating the reduced bases at the nodes of the reference mesh, the shape functions play a crucial role. We found that using elements with higher-order shape functions has a beneficial effect on the accuracy of the interpolated basis. We therefore recommend using at least elements with quadratic shape functions. With regard to mesh morphing by using the spring analogy with elastic hardening, we follow the idea of \cite{Nazem2021} when using triangular elements with quadratic shape functions: The element is divided into four sub-elements, which are treated as four elements for the spring analogy. 

\section{Results}
\label{sec:Results}

The proposed methods for pMOR by matrix interpolation for varying underlying meshes are tested on three structures: a beam-shaped plate, a plate with a circular hole, and a plate with an elliptic hole. In all cases, triangular plane stress elements with quadratic shape functions are used. For the plate, a structured mesh was created manually; for the other two examples, the toolbox \verb+distmesh+ \cite{Persson2004} was used for mesh generation. This toolbox only creates triangular elements with linear shape functions; thus, nodes at the mid-points of all edges of the triangles were added after the automatic mesh generation to obtain triangular elements for quadratic shape functions. MATLAB 2023a was used on an Intel Xeon E5-2660v3 processor with 2.6 GHz for all computations. For the problems with one-dimensional parameter spaces, the reduced operators were interpolated using spline interpolation; for the two-dimensional parameter space, ridge regression with polynomials with all terms up to order 3 in each variable and a ridge parameter of $\lambda=10^{-5}$ was used, which proved to be good choices. \\

\subsection{Beam-shaped Plate}
\label{sec:Plate}

\begin{figure}
    \centering
    \includegraphics[width=0.9\linewidth]{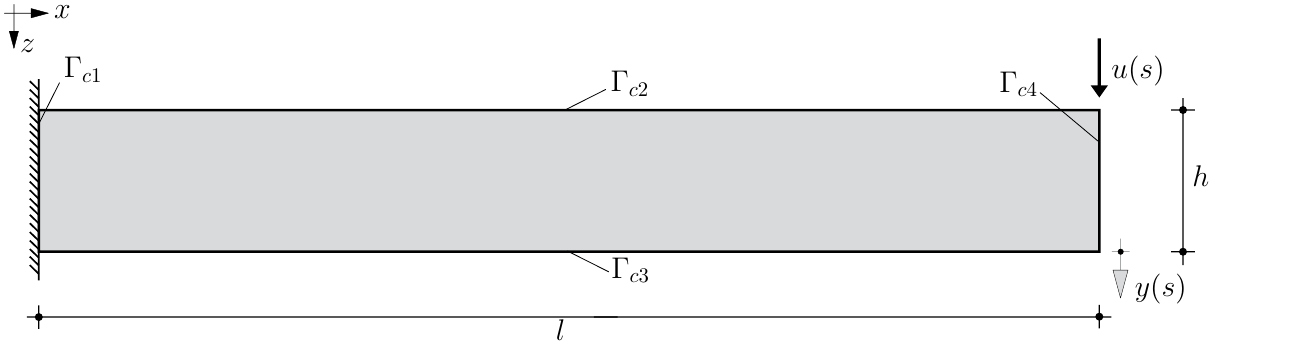}
    \caption{Beam-shaped plate.}
    \label{fig:Plate}
\end{figure}

We first investigate the framework on an academic example, namely a beam-shaped plate depicted in  Fig. \ref{fig:Plate}. This example is based on the benchmark problem presented in \cite{Panzer2009}. As parameter, the length of the beam is used and varied between $0.8 \, \mathrm{m}$ and $1.2 \, \mathrm{m}$. The remaining parameters of the model are fixed and chosen as shown in Table~\ref{tab:ParametersPlate}. All degrees of freedom at the left edge of the plate are fixed. On the top right corner of the plate, a time-harmonic force is applied. The excitation frequency is varied in the interval [1, 5000] Hz in steps of 1 Hz. The output is the vertical displacement of the plate at the bottom right corner. The plate is meshed with quadratic triangular elements, with a mesh size of 0.02, resulting in a size of the full-order model ranging from 1782 to 2662 degrees of freedom, depending on the length.

\begin{table}[htb]
	\caption{Geometry and material parameters of the beam-shaped plate and hyperparameter values used for the adaptive sampling.}\label{tab:ParametersPlate}
	\centering
	\begin{tabular}{ c c c c c c } 
		Parameter & Range/Value & Unit & & Hyperparameter & Value \\  \cmidrule{1-3} \cmidrule{5-6}
		Length $l$ & $[0.8, \, 1.2]$ & $\mathrm{m}$ & & $\theta_{lT}$ & $10^\circ$ \\  
		Height $h$ & $0.1$ & $\mathrm{m}$ & & $\theta_{uT}$ & $85^\circ$ \\
		Thickness $t$ & $0.01$ & $\mathrm{m}$ & & $d_{lT}$ & $0.1$ \\
		Young's modulus $E$ & $2.1\cdot 10^{11}$ & $\mathrm{N/m^2}$ & & $d_{uT}$ & $0.2$ \\
		Poisson's ratio $\nu$  & $0.3$ & - & & $d_N$ & 0 \\
		Density $\rho$ & $7860$ & $\mathrm{kg/m^3}$ & & minimum amount of & \multirow{2}{*}{4} \\
		Rayleigh damping $\alpha$ & $8$ & $\mathrm{1/s}$ & & samples per cluster  & \\
		Rayleigh damping $\beta$ & $8\cdot 10^{-6}$ & $\mathrm{s}$ & & &
	\end{tabular}
\end{table}

We employ the adaptive sampling and clustering algorithm proposed in \cite{Resch-Schopper2024} with the extension for varying underlying meshes to generate a consistent parametric reduced-order model (pROM). For the adaptive sampling algorithm, the reduced operators are computed for two initial samples at $l = 0.8\, \mathrm{m}$ and $l = 1.2 \, \mathrm{m}$ using modal truncation. Based on the rule of thumb to use all eigenmodes up to twice the regarded frequency range, the reduced size is set to 16. The hyperparameters of the adaptive sampling algorithm are set as shown in Table~\ref{tab:ParametersPlate}. In the adaptive sampling and for transforming the reduced operators after the sampling, morphing by the spring analogy with elastic hardening and radial basis function morphing are used. As characteristic points, the four edges of the plate are chosen, so that the following displacements are prescribed for morphing:
\begin{align}
    \bar{\q}_x &= 0 \quad \text{on } \Gamma_{c1}, \\
    \bar{\q}_z &= 0 \quad \text{on } \Gamma_{c1}, \, \Gamma_{c2} \text{ and } \Gamma_{c3}, \\
    \bar{\q}_x &= \Delta l \quad \text{on } \Gamma_{c4}.
\end{align}
Here, $\Delta l$ is the difference between the lengths of the two underlying geometries. For morphing by the spring analogy, the prescribed displacement is imposed in ten steps with a tenth of the total magnitude in each step, which proved to be sufficiently many steps to avoid collapsing elements. In radial basis function morphing, the linear kernel function $\psi(\|\x - \x_i \|_2) = \| \x - \x_i \|_2$ is used. \\
For comparison, two alternative pROMs are computed: once using the basis transformation proposed in \cite{Amsallem2016} and once by adding zero matrices to the sampled reduced bases as suggested in \cite{Roy2021}, so that the steps of pMOR by matrix interpolation outlined in Eqs. (\ref{eq:Vg}) - (\ref{eq:TransfOp}) can be performed. For both, the samples obtained from the adaptive sampling algorithm with radial basis function morphing are used. The accuracy of the pROMs is assessed for 1001 linearly spaced test points in the range from $0.8 \, \mathrm{m}$ to $1.2 \, \mathrm{m}$ by computing the mean of the relative error between the full-order output $y(s)$ and the predicted output $y_r(s)$:
\begin{equation}
    \varepsilon = \frac{1}{n_s} \sum_{i=1}^{n_s} \frac{|y(s_i) - y_r(s_i)|}{|y(s_i)|}.
\end{equation}
Here, $s_i$ are the discrete frequencies for which the frequency response function is evaluated and $n_s$ is the number of discrete frequencies used.
\begin{figure}[htb]
	\centering
	\includegraphics[width=\textwidth]{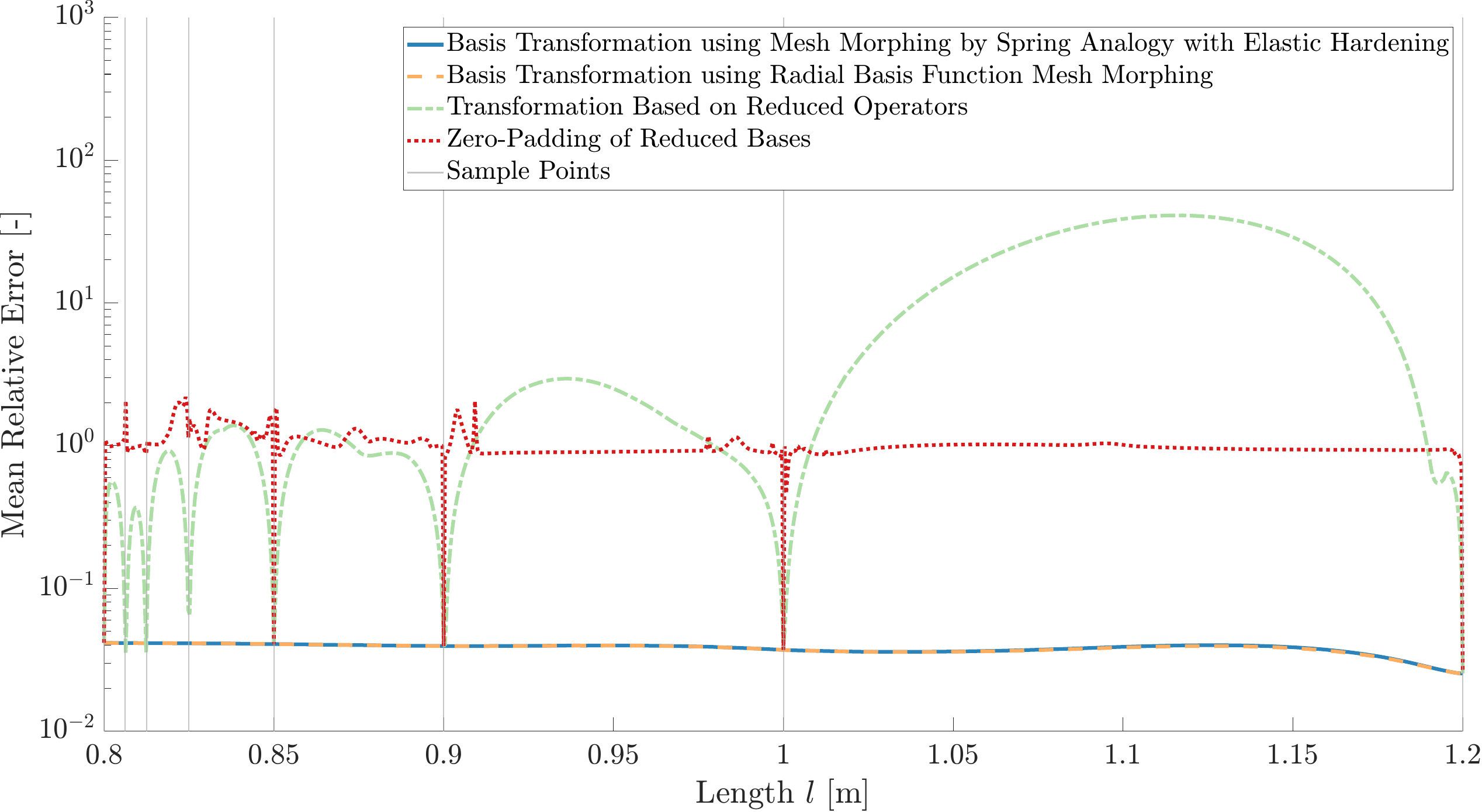}
	\caption{Mean relative error of the predicted reduced order models for a beam-shaped plate.}
	\label{fig:MRE_plate}
\end{figure}

\begin{figure}[htb]
	\centering
	\includegraphics[width=\textwidth]{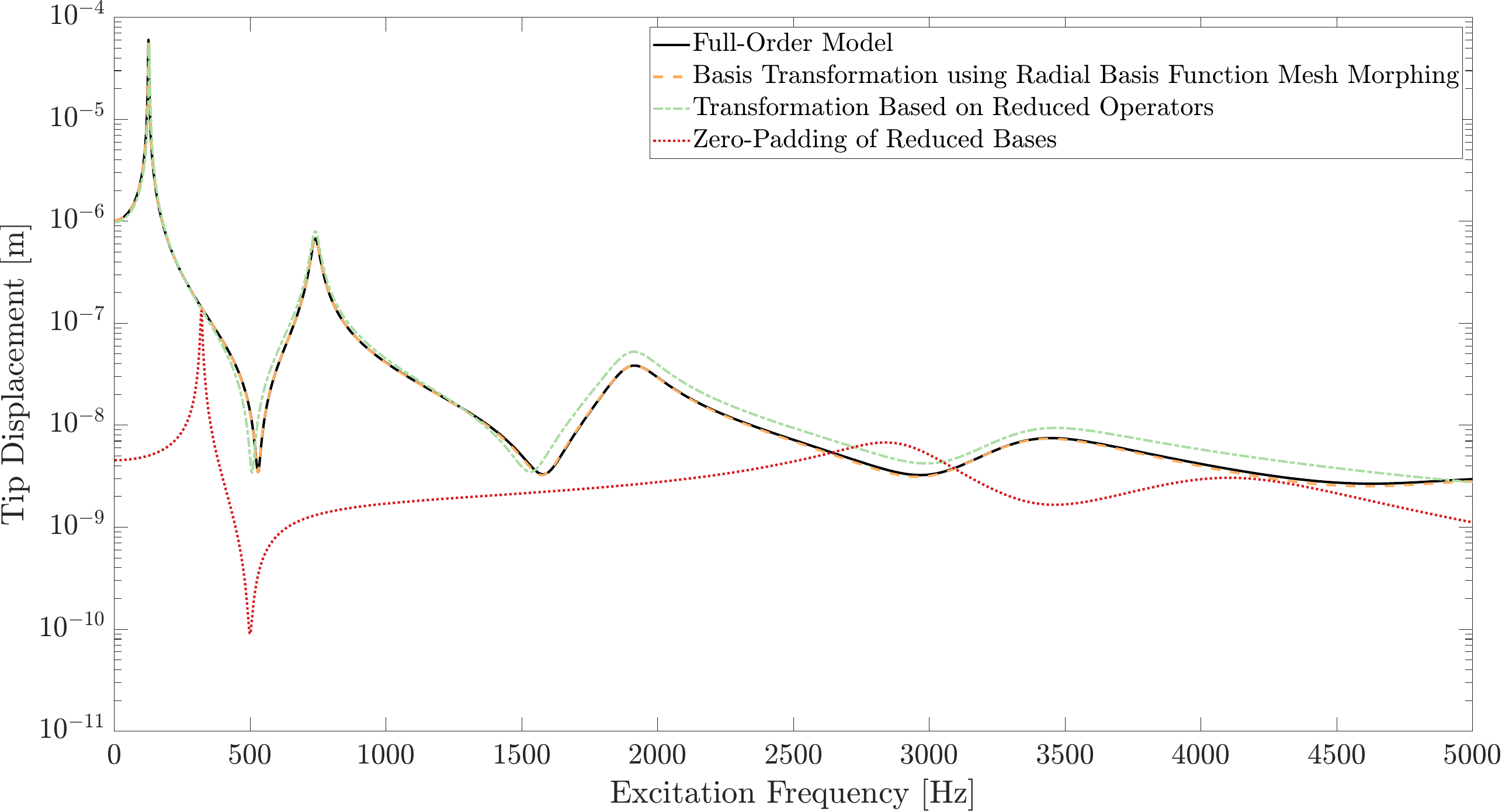}
	\caption{Frequency response function of the full-order model and the predicted reduced-order models for the beam-shaped plate.}
	\label{fig:FRF_Plate}
\end{figure}

Fig. \ref{fig:MRE_plate} shows the mean relative error over the parameter range for the four pROMs. For this small reduced size, all eigenmodes can be matched for the whole parameter range. Therefore, the adaptive sampling algorithm computes in total 8 samples to ensure small subspace angles between the reduced bases of neighboring samples, but does not divide the parameter space into multiple regions. It can be seen that the accuracy of the two proposed pROMs, for which the basis transformation is performed using mesh morphing, is almost identical. In general, the accuracy of these two pROMs does not deviate from the accuracy at the sample points, indicating that the reduced operators are consistent and can be interpolated accurately. \\
The error of the pROM based on the transformation of the reduced operators increases significantly between the sample points, indicating that the resulting parametric dependency of the transformed reduced operators cannot be interpolated well with polynomials. The reason for this behavior is that the transformation based on the reduced operators cannot distinguish between a change in the full operators and a change in the reduced basis and tries to minimize both simultaneously. This can lead to a parametric dependency that is difficult to interpolate. \\
The pROM, for which the transformation has been performed by zero-padding the reduced bases, shows high errors throughout the whole parameter space. When adding zeros to the reduced bases, the displacement field they represent is altered in comparison to other reduced bases, which prevents a meaningful transformation to a common coordinate system. Consequently, a physical interpolation of the reduced operators is not possible, which leads to these high errors. \\
For a more detailed analysis, the frequency response function (FRF) of the full model and the pROMs is shown in Fig. \ref{fig:FRF_Plate} for a length of $l = 0.81 \, \mathrm{m}$, which is between the second and the third sample point. Since the predicted FRFs of the two pROMs based on mesh morphing are almost identical, only the result for radial basis function morphing is shown. As already indicated by the mean relative error, this pROM accurately captures the behavior of the FOM. For the pROM resulting from the transformation based on the reduced operators, the FRF deviates slightly because the reduced operators cannot be interpolated accurately. For the pROM based on zero-padding, a meaningful transformation and interpolation of the reduced operators is not possible, which is why the predicted FRF deviates significantly from the full model's behavior.

\subsection{Plate With Hole}
\label{sec:PlateHole}

\subsubsection{One-dimensional Parameter Space}

\begin{figure}
    \centering
    \includegraphics[width=0.65\linewidth]{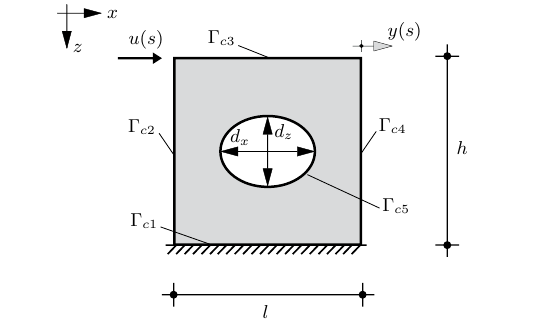}
    \caption{Plate with hole.}
    \label{fig:PlateHole}
\end{figure}

Next, we investigate the framework on a plate with a hole depicted in Fig. \ref{fig:PlateHole}. We first use the same diameter in $x$- and $z$-direction as parameter, such that $d_x = d_z = d$. The remaining parameters of the model are fixed and chosen as shown in Table~\ref{tab:ParametersPlateHole}. All degrees of freedom at the bottom edge of the plate are fixed. On the top left corner of the plate, a time-harmonic excitation is applied; the output is the horizontal displacement of the plate at the top right corner. The plate is meshed with quadratic triangular elements with an element size of $0.02\, \mathrm{m}$. For a finer mesh at the circle, the following refinement function is used for the relative element size distribution:
\begin{equation}
    h(x,z) = 0.05 + 0.3 \cdot \sqrt{(x - x_c)^2 + (z - z_c)^2} - \frac{d}{2}.
\end{equation}
Here, $x_c$ and $z_c$ are the $x$- and $z$-coordinate of the center of the circular hole. For more information about the definition of the element size and refinement in \verb+distmesh+, the reader is referred to \cite{Persson2004}. Depending on the diameter, the resulting mesh has between 4126 and 5660 degrees of freedom. 

\begin{table}[htb]
	\caption{Geometry and material parameters of the plate with hole and hyperparameter values used for the adaptive sampling.}\label{tab:ParametersPlateHole}
	\centering
	\begin{tabular}{ c c c c c c } 
		Parameter & Range/Value & Unit & & Hyperparameter & Value \\  \cmidrule{1-3} \cmidrule{5-6}
		Diameter $d$ & $[0.2, \, 0.6]$ & $\mathrm{m}$ & & $\theta_{lT}$ & $10^\circ$ \\  
		Height $h$ & $1$ & $\mathrm{m}$ & & $\theta_{uT}$ & $85^\circ$ \\
		Thickness $t$ & $0.01$ & $\mathrm{m}$ & & $d_{lT}$ & $0.05$ \\
		Young's modulus $E$ & $2.1\cdot 10^{11}$ & $\mathrm{N/m^2}$ & & $d_{uT}$ & $0.2$ \\
		Poisson's ratio $\nu$  & $0.3$ & - & & $d_N$ & 0 \\
		Density $\rho$ & $7860$ & $\mathrm{kg/m^3}$ & & minimum amount of & \multirow{2}{*}{4} \\
		Rayleigh damping $\alpha$ & $8$ & $\mathrm{1/s}$ & & samples per cluster  & \\
		Rayleigh damping $\beta$ & $8\cdot 10^{-6}$ & $\mathrm{s}$ & & &
	\end{tabular}
\end{table}

For the adaptive sampling algorithm, the reduced operators are computed for two initial samples at $d = 0.2\, \mathrm{m}$ and $d = 0.6 \, \mathrm{m}$ using modal truncation with a reduced size of 50. This again refers to the number of eigenmodes within twice the frequency range of interest. The hyperparameters of the adaptive sampling algorithm are set as shown in Table~\ref{tab:ParametersPlateHole}. Again, morphing by the spring analogy with elastic hardening and radial basis function morphing are used with the same settings as for the previous example and compared to the pROMs built by the approaches proposed in \cite{Amsallem2016} and \cite{Roy2021}. As characteristic points, the circumference of the hole is chosen in addition to the four edges of the plate. This results in the following displacements to be prescribed for morphing:
\begin{align}
    \bar{\q}_z &= 0 \quad \text{on } \Gamma_{c1}, \text{ and } \Gamma_{c3},\\
    \bar{\q}_x &= 0 \quad \text{on } \Gamma_{c2}, \text{ and } \Gamma_{c4}, \\
    \bar{\q} &= \frac{\Delta d}{2} \cdot (\mathbf{x} - \mathbf{x}_c)  \quad \text{on } \Gamma_{c5}.
\end{align}
Here, $\Delta d$ is the difference in the diameters of the two underlying geometries, $\mathbf{x} = [x, z]^\intercal$ is a point on $\Gamma_{c5}$, and $\mathbf{x}_c = [x_c, z_c]^\intercal$ corresponds to the center of the hole. \\
The accuracy of the pROMs is assessed for 1001 linearly spaced sample points in the range from $0.2 \, \mathrm{m}$ to $0.6 \, \mathrm{m}$ by computing the mean of the relative error between the full-order output $y(s)$ and the predicted output $y_r(s)$ of the pROMs for 5000 linearly spaced values in the frequency range $[1, \, 5000]$ Hz. 
\begin{figure}[htb]
	\centering
	\includegraphics[width=\textwidth]{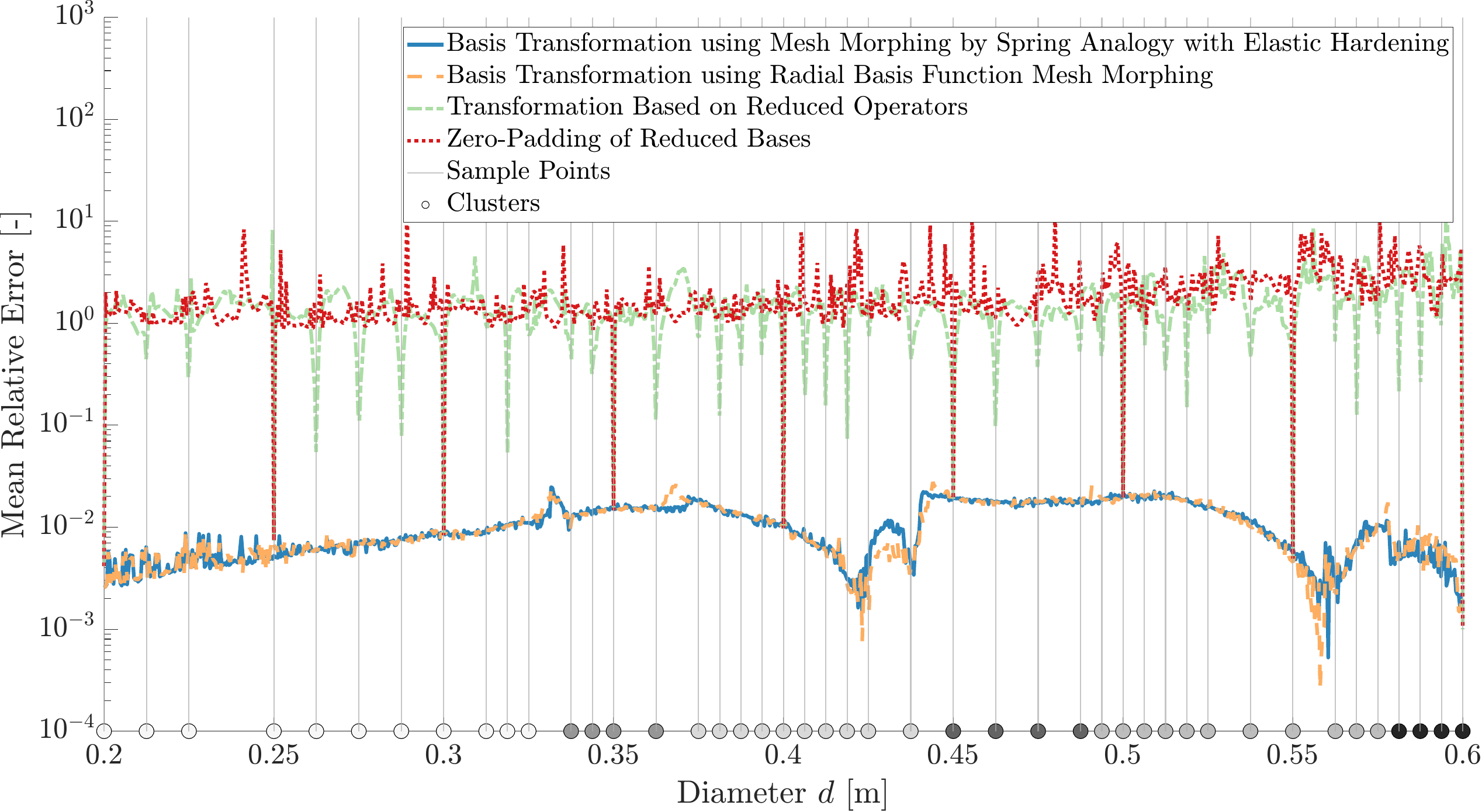}
	\caption{Mean relative error of the predicted reduced order models for the plate with a circular hole.}
	\label{fig:MRE_PlateCircularHole}
\end{figure}

\begin{figure}[htb]
	\centering
	\includegraphics[width=\textwidth]{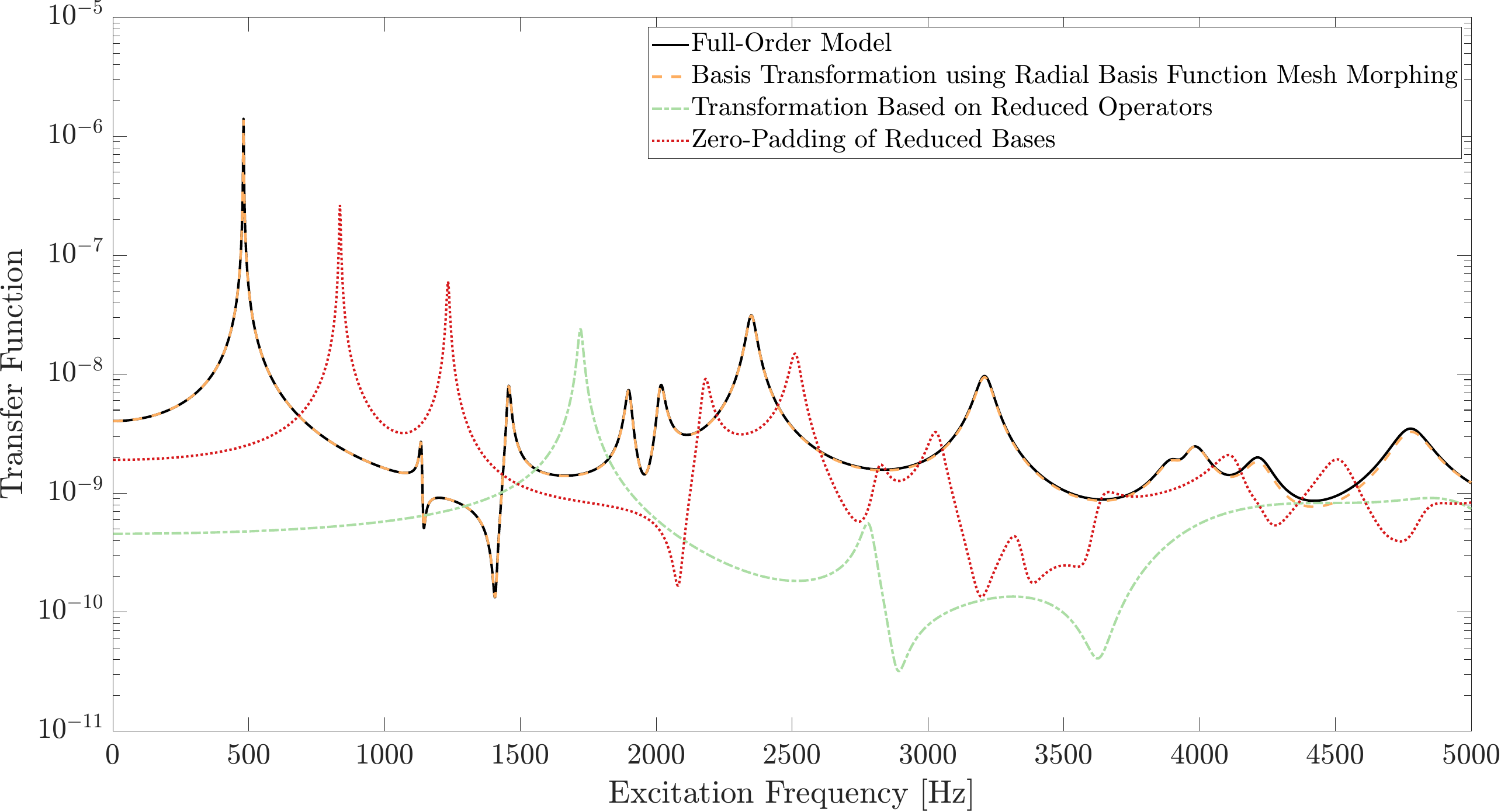}
	\caption{Frequency response function of the full-order model and the predicted reduced-order models for the plate with a circular hole.}
	\label{fig:FRF_PlateCircularHole}
\end{figure}

In Fig. \ref{fig:MRE_PlateCircularHole}, the mean relative error of the four pROMs is shown over the parameter space. With radial basis function morphing, the adaptive sampling algorithm computes in total 44 samples. This time, the sampled reduced bases cannot be transformed to the same common coordinate system for all samples. Therefore, the algorithm clusters the samples into 6 groups, which are shown on the $x$-axis of Fig. \ref{fig:MRE_PlateCircularHole}. Using morphing by the spring analogy with elastic hardening leads to 47 samples with a similar distribution and clustering. Consequently, the results of both pROMs based on morphing are again very similar and show a consistently high level of accuracy across the entire parameter space. This time, however, the accuracy between the sample points deviates slightly from the accuracy at the sample points. This may be explained by two reasons. Firstly, the mesh morphing methods used are not exact. This means that a morphed node might be moved to a position that is not equivalent to its relative position in the initial mesh, which introduces further inconsistencies in the reduced operators. However, as shown in Fig. \ref{fig:MRE_PlateCircularHole}, it only marginally affects the accuracy of the pROMs. In the previous example shown in Section~\ref{sec:Plate}, this did not occur because the exact morphing consists of a linear scaling of the $x$-coordinate, which is well represented by both the spring analogy and the radial basis function morphing. Secondly, \verb+distmesh+ does not always generate the exact same mesh for the same geometry. Therefore, slight differences in the results can be introduced due to the different meshes. \\
The two alternative approaches from the literature for generating a pROM exhibit high errors across the entire parameter space. The reason for this behavior is that the reduced bases of all samples cannot be transformed to the same coordinate system in this case. This, however, cannot be detected by these two approaches. The reduced operators then contain inconsistencies, and interpolating them lacks physical meaning, which leads to these high errors. \\
This behavior can also been seen in detail in Fig.~\ref{fig:FRF_PlateCircularHole}, which shows the FRFs of the pROMs for $d = 0.33 \, \mathrm{m}$. This refers to the point between the first and the second cluster, where the two pROMs based on mesh morphing show the highest error. The mean relative error for this point is still below 2.5 \% though, which can also be seen in the accurate representation of the full model's behvaior. As already indiciated by the mean relative error, the two other pROMs cannot capture the full model's behavior and their FRFs deviate significantly.

\subsubsection{Two-dimensional Parameter Space}

We now use the diameters $d_x$ and $d_z$ of the plate shown in Fig. \ref{fig:PlateHole} as two separate parameters and vary them both in the range $[0.2, \, 0.4]$ m. The other parameters of the plate remain the same as listed in Tab. \ref{tab:ParametersPlateHole}. We again use quadratic triangular elements with an element size of 0.02 for meshing, but now use the following relative element size distribution function for a refinement at the ellipse:
\begin{equation}
      h(x,z) = 1 + 0.5 \cdot \sqrt{\frac{(x - x_c)^2}{(0.5 \cdot d_x)^2} + \frac{(z - z_c)^2}{(0.5 \cdot d_z)^2}} - 1.
\end{equation}
The resulting meshes comprise between 4486 and 9190 degrees of freedom, depending on the diameters. As initial samples, the four corner points of the parameter space are used, and reduced with modal truncation using a reduced size of 50. The hyperparameters for the adaptive sampling and clustering algorithm stay mostly the same as listed in Table~\ref{tab:ParametersPlateHole}. Only $d_{lT}$ is set to 0.1 so that the border between separate consistent regions is sampled a bit coarser and the minimum amount of samples per cluster is set to 1 because there is no minimum required amount of samples for ridge regression. Furthermore, finding the border of the regions (Algorithm 4.5 in \cite{Resch-Schopper2024}) is skipped here in order to limit unnecessarily many samples at the border between neighboring consistent regions. For the morphing methods, the following displacements are prescribed:
\begin{align}
    \bar{\q}_z &= 0 \quad \text{on } \Gamma_{c1}, \text{ and } \Gamma_{c3},\\
    \bar{\q}_x &= 0 \quad \text{on } \Gamma_{c2}, \text{ and } \Gamma_{c4}, \\
    \bar{\q}_x &= \frac{\Delta d_x}{2} \cdot (x - x_c)  \quad \text{on } \Gamma_{c5}, \\
    \bar{\q}_z &= \frac{\Delta d_z}{2} \cdot (z - z_c)  \quad \text{on } \Gamma_{c5}.
\end{align}
Here, $\Delta d_x$ and $\Delta d_z$ are the differences in diameters in the $x$- and $z$-direction, respectively. \\
The accuracy of the pROMs is assessed for $51 \times 51$ linearly spaced sample points in the range from 0.2~m to 0.4~m by computing the mean of the relative error between the full-order output $y(s)$ and the predicted output $y_r(s)$ of the pROMs for 5000 linearly spaced values in the frequency range [1, 5000] Hz.
\begin{figure*}[htb]
	\centering
	\begin{subfigure}[b]{0.47\textwidth}
		\centering
		\includegraphics[width=\textwidth]{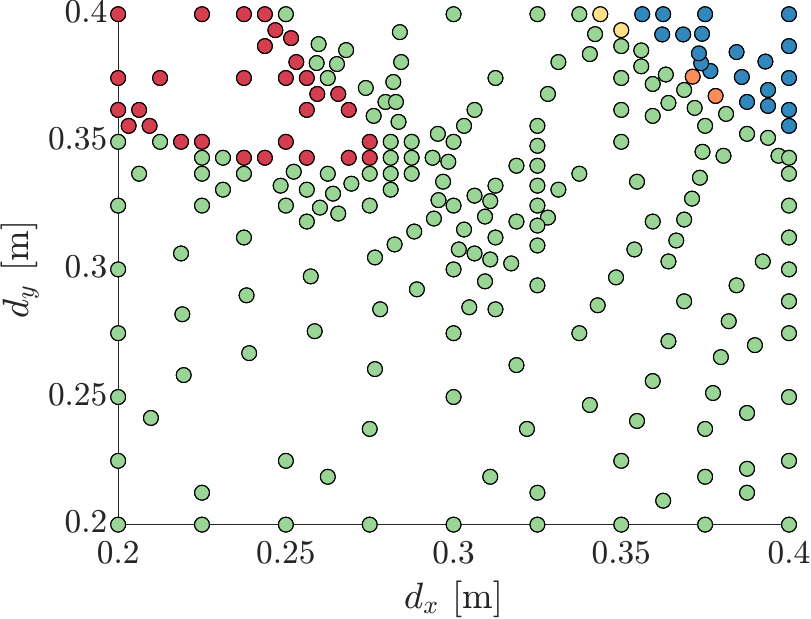}
		\caption{Final distribution and clustering of samples.}  
		\label{fig:Samples_PlateEllipticHole_SAEH}
	\end{subfigure}
	\hfill
	\begin{subfigure}[b]{0.51\textwidth}  
		\centering 
		\includegraphics[width=\textwidth]{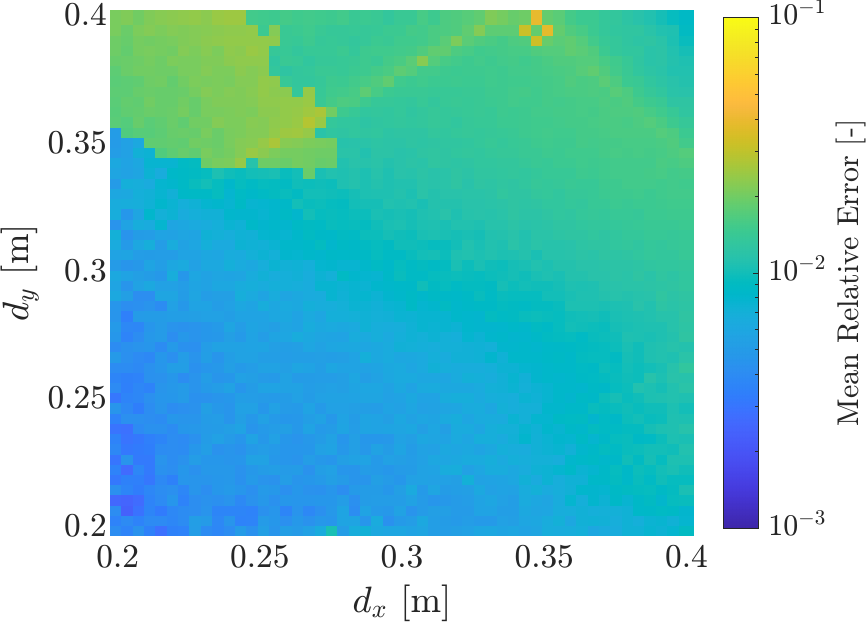}
		\caption{mean relative error of the parametric reduced model.}  
		\label{fig:MRE_PlateEllipticHole_SAEH}
	\end{subfigure}
	\caption{Results of the adaptive sampling and predicted reduced models for the plate with the elliptic hole using the spring analogy with elastic hardening for morphing.}  
	\label{fig:MRE_PlateEllipticHole}
\end{figure*}

Fig. \ref{fig:MRE_PlateEllipticHole} shows the results of the adaptive sampling and the predicted reduced models for the pROM based on morphing by the spring analogy with elastic hardening. The sampling algorithm placed 233 samples in addition to the initial samples and clustered them into 5 groups. In the whole parameter space, the error of the predicted reduced-order models is below 4 \%, in most regions even well below this. There are different levels of accuracy in the different regions, and to a degree, generally higher errors at the border between separate regions. Using radial basis function morphing yields very similar results, both in the distribution of the 247 added samples clustered into 6 groups and in the mean relative error. The pROMs based on the transformation of the reduced operators, respectively, on an adjustment of the reduced bases by zero-padding, suffer from the inconsistencies in the reduced operators, similarly to the previous example, and show errors of around 100~\% in the entire parameter space.

\subsection{Comments on the computational cost}
\label{sec:Runtime}

\begin{figure}[htb]
	\centering
	\includegraphics[width=0.5\textwidth]{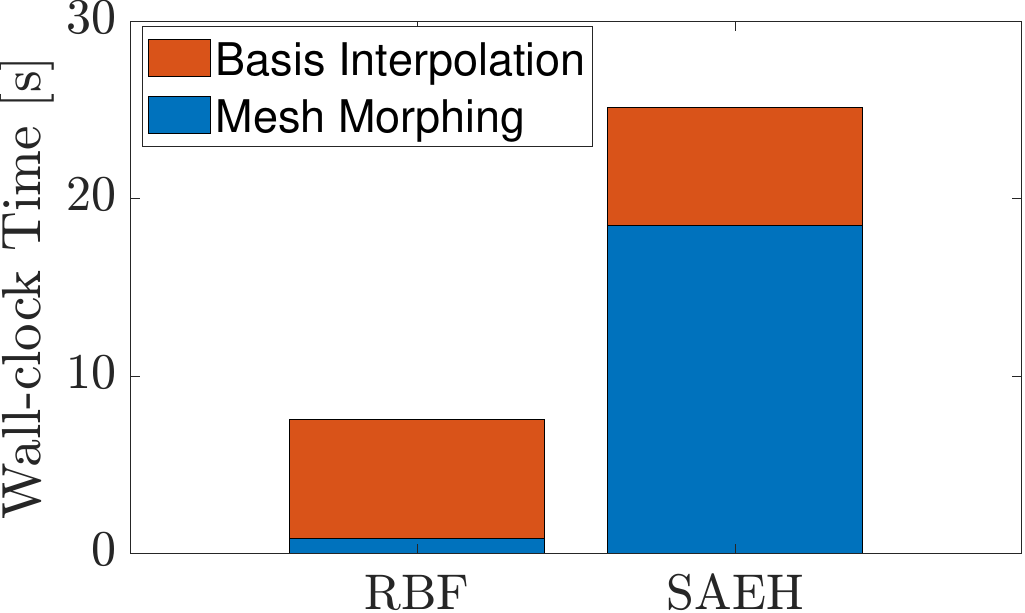}
	\caption{Required time to perform the basis interpolation and the mesh morphing with radial basis function (RBF) morphing and morphing by spring analogy with elastic hardening (SAEH), respectively, for the plate with an elliptic hole.}
	\label{fig:Runtime}
\end{figure}

Advantageously, all costly computations can be performed in the offline phase; in the online phase, the reduced-order models are directly predicted and can be solved efficiently. However, in the offline phase, morphing the meshes and interpolating the reduced bases to a morphed mesh must be performed multiple times. During the adaptive sampling, this is required every time the subspace angles for two samples needs to be computed, for which it has not been computed before. Furthermore, to perform the transformation of the reduced operators within a consistent region after the adaptive sampling algorithm, mesh morphing and basis interpolation must be applied to all samples of that cluster, except for the reference sample. Therefore, the runtime of these operations shall be commented on here. \\
As the reference sample, $d_x = d_z = 0.4 \, \mathrm{m}$ is chosen, which leads to the densest mesh of all samples in the parameter space with more than 9000 degrees of freedom. For performing the morphing and basis interpolation, a reduced basis of size $r=50$ is computed using modal truncation for the configuration $d_x = 0.38, \, \mathrm{m}, \, d_z = 0.34 \, \mathrm{m}$, which shall be expressed in terms of the reference mesh. The computations are performed using MATLAB 2023a on an Intel Xeon E5-2660v3 processor with 2.6 GHz, utilizing parallel computing with 20 cores where possible. \\
In Fig. \ref{fig:Runtime}, the wall-clock time required to perform mesh morphing and basis interpolation is shown. It can be seen that mesh morphing using radial basis functions can be performed in less than one second, which is significantly faster than morphing based on the spring analogy. This is because for fitting the radial basis function model, only Eq.~\eqref{eq:RBF_ES} must be solved. The size of the matrix $\mathbf{A}$ is the number of characteristic nodes, which is slightly greater than 200 for the above case, and the number of columns of $\mathbf{B}$ is the dimension of the structure plus one, so three in this example. Therefore, solving this equation system is computationally inexpensive and must be done only once per morphing. Predicting the position of the non-characteristic nodes can then be performed efficiently in parallel. For performing mesh morphing using the spring analogy by elastic hardening, the prescribed displacement must be applied in multiple steps. Each time the mesh is deformed, the new spring stiffnesses must be computed and assembled. For the investigated example, this takes around 18 seconds when the prescribed displacement is applied in ten steps. Thereby, most of the time is required to solve the equation system to compute the intermediate deformed mesh. The approach for interpolating the reduced basis is identical for both morphing approaches and takes around 6 seconds for the investigated example for 50 basis vectors. Thereby, most of the time is required to evaluate the shape functions for all nodes of the morphed mesh. This step, as well as finding the element in which the nodes of the morphed mesh lie, can be performed in parallel, which could also decrease the runtime on a machine with more cores.

\section{Conclusion}
\label{sec:conclusion}

A method for generating consistent parametric reduced-order models by matrix interpolation has been presented for cases where the mesh of the underlying structure varies. The key idea is to understand the reduced basis vectors as continuous displacement fields that can be represented in terms of different discretizations. By mesh morphing and basis interpolation, the reduced bases of samples with varying meshes can be described in terms of one reference mesh. This allows a comparison of the sampled bases, for example through subspace angle analysis, which is important to detect strong changes in the reduced bases and minimize them. If not treated, these strong changes lead to inconsistencies in the reduced operators, which prevent a meaningful and accurate interpolation of the entries. \\
For the mesh morphing, two strategies were implemented and tested, namely morphing by spring analogy with elastic hardening and radial basis function morphing. Both approaches showed similarly high accuracies for one- and two-dimensional test cases of a beam-shaped plate and a plate with a hole. They performed significantly better than two existing approaches for pMOR by matrix interpolation for varying underlying meshes, including transformations based on the reduced operators or adding zero matrices to the reduced bases to embed them in the same high dimension. If inconsistencies appear in the training data due to large parameter variations, these two existing approaches proved to be unsuitable. Regarding the computational effort, radial basis function morphing can be performed in a fraction of the time required for morphing by the spring analogy and thus leads to significantly smaller offline costs. 
These results demonstrate that the proposed framework provides a robust and computationally efficient foundation for consistent parametric model reduction with non-matching meshes, enabling applications to more complex, multi-parameter, and three-dimensional problems.
 
\section*{Acknowledgments}
R. Rumpler gratefully acknowledges the financial support provided by the Swedish Research Council (VR Grant 2021-05791).

\section*{Conflict of interest}
The authors declare no potential conflict of interest.

\section*{Data Availability}
The code for the adaptive sampling and clustering algorithm for matching underlying meshes can be accessed via the following github repository: \url{https://github.com/SebastianResch-Schopper/Consistent-Parametric-Model-Order-Reduction-by-Matrix-Interpolation.git} \\
The data for varying underlying meshes will be made available upon request.

\section*{CRediT authorship contribution statement}
\textbf{Sebastian Resch-Schopper:}  Writing - original draft, Methodology, Software, Validation, Visualization. \textbf{Romain Rumpler:} Supervision, Writing - review \& editing. \textbf{Gerhard Müller:} Supervision, Writing - review \& editing.

\bibliography{bibtex/references} 

\end{document}